\documentclass [11pt]{amsart}
\usepackage {amsmath, amssymb, amscd}

\newtheorem{theorem}{Theorem}
\newtheorem {definition}{Definition}
\newtheorem {lemma}{Lemma}
\newtheorem {proposition}{Proposition}

\newtheorem {conjecture}{Conjecture}
\newtheorem {remark}{Remark}

\def\go {{\mathcal{G}_0}}
\def\gau {{\mathcal{G}}}
\def\ss {{S}}
\def\zz {{\mathbb{Z}}}
\def\rr {{\mathbb{R}}}
\def\cc {{\mathbb{C}}}
   
\def\spc {{\mathfrak{s}}}

\def\dif {{d}}

\def\mi {{i}}
\def\inte {{\text{int}}}
\def\spa {{\mi\Omega^1(Y) \oplus \Gamma(W_0)}}
\def\swf {{\text{SWF}}}
\def\swfn {{\text{\ttfamily swf}}}
\def\dir {{\not\!\partial}}

\def\vp {{\varphi}}
\def\del {{\partial}}
\def\delt {{\frac{\partial}{\partial t}}}

\def\delp {{(\del/\del t)}}
\def\os {{\text{Ozsv\'ath-Szab\'o}}}
\def\fin {{\hspace{15pt} \square}}

\def\pr {{\text{pr}}}

\def\colim {{\text{colim}}}
\def\im {{\hspace{2pt} \text{Im }}}

\def\sp {{\mathfrak{S}}}
\def\pros {{\text{Pro-}\mathfrak{S}'}}

\def\from {{\leftarrow}}
\begin{document}

\title [Periodic Floer pro-spectra] {Periodic Floer pro-spectra from the 
Seiberg-Witten equations}
\author [Peter B. Kronheimer and Ciprian Manolescu]{Peter B. 
Kronheimer and Ciprian Manolescu}
\date {September 24, 2003}

\begin {abstract} Given a three-manifold with $b_1 =1$ and a nontorsion 
$spin^c$ structure, we use finite dimensional approximation to construct
from the Seiberg-Witten equations two invariants in the form of a periodic
pro-spectra. Various functors applied to these invariants give different 
flavors of Seiberg-Witten Floer homology. We also construct 
stable homotopy versions of the relative Seiberg-Witten invariants 
for certain four-manifolds with boundary.
\end {abstract}

\address {Department of Mathematics, Harvard University\\ 1 Oxford 
Street, Cambridge, MA 02138}
\email {kronheim@math.harvard.edu}

\address {Department of Mathematics, UCLA\\ 520 Portola Plaza, Los Angeles, CA 90095}
\email {cm@math.ucla.edu}

\maketitle

\section {Introduction}

In \cite{CJS}, Cohen, Jones, and Segal posed the question of constructing
a ``Floer homotopy type.'' They conjectured that Floer homology (in either  
of the two variants known at the time, symplectic or instanton) should be 
the homology of an object called a \emph {pro-spectrum.} However, the 
passage from homology to homotopy in either of these cases seems a  
difficult task at the moment. Because of their remarkable compactness 
properties, the Seiberg-Witten equations are better suited for this task. Using 
the technique of finite dimensional approximation, Furuta and Bauer were 
able to define stable homotopy invariants for four-manifolds (\cite{B}, 
\cite{Fu1}, \cite{Fu2}). In \cite{M}, the second author has associated to 
each rational homology 3-sphere a certain equivariant spectrum whose homology 
is the Seiberg-Witten Floer homology.

This paper is a continuation of \cite{M}. Here we define the stable
homotopy generalization of the Seiberg-Witten Floer homology for 
closed, oriented 3-manifolds $Y$ with $b_1(Y)=1,$ together with $spin^c$ 
structures $\spc.$ We assume that $c_1(\spc),$ the first Chern class of 
the corresponding determinant line bundle, is nontorsion. 

It turns out that the usual homotopy category $\sp$ of spectra 
is not good enough to support our invariants. To ensure the existence 
of good colimits, we divide out the sets of morphisms in $\sp$ by the 
class of phantom maps, and we call the resulting category $\sp'.$ 
Furthermore, to ensure the existence of good inverse limits, we enlarge 
the class of objects to include \emph{pro-spectra}. The resulting category 
is called $\pros.$ Its exact definition will be given in section~\ref{sec:g}.

The main result of this article is: 
\begin {theorem} \label {nontorsion}
a) Given a Riemannian metric $g$ on $Y$ and a $spin^c$ connection 
$A_0,$ finite dimensional approximation
for the Seiberg-Witten map produces an invariant $\swf(Y, \spc, g, A_0)$ 
in the form of a pointed $S^1$-equivariant pro-spectrum well-defined up to
canonical equivalence in $\pros.$ 

When the metric $g$ or the connection $A_0$ change, the invariant 
$\swf(Y, \spc, g, A_0)$ can change only by suspending or desuspending by a 
finite dimensional complex representation of $S^1.$ 

Furthermore, the pro-spectrum is periodic
modulo $\ell,$ where
$$ \ell = g.c.d. \{ (h \cup c_1(\spc))[Y] \hspace{3pt} | h \in H^1(Y;
\zz)\},$$ in the sense that there is a complex representation $E$ of $S^1$ 
of real dimension $\ell$ and a natural equivalence:
$$\swf \to \Sigma^{E} (\swf).$$

b) There is a different version $\swf_0(Y, \spc, g, A_0)$ which is a 
spectrum obtained by 
doing finite dimensional approximation for the Seiberg-Witten map with a 
nonexact perturbation in the cohomology class $-c_1(\spc).$ 
It has the same properties as the ones 
mentioned above for $\swf(Y, \spc, g, A_0),$ 
and there is a natural morphism:
$$ j: \swf_0 \to \swf $$
which induces isomorphisms on homology and on the equivariant Borel 
homology.

c) Under change of orientation, $\swf(Y, \spc, g, A_0)$ and $\swf(\bar Y, 
\bar \spc, g, A_0)$ are duals, but the analogous statement for $\swf_0$ is 
false.
\end {theorem}

Our construction of $\swf$ runs parallel to the one for 
rational homology spheres in \cite{M}. However, there is an 
important complication, given by the loss of compactness 
of the Seiberg-Witten moduli space. If we work in Coulomb gauge, the 
solutions of the Seiberg-Witten equations are the critical points of the 
Chern-Simons-Dirac functional 
$$ CSD : V=(\Omega^1(Y; i\rr)/ \im d) \oplus \Gamma (W_0) \to \rr, $$
where $W_0$ is a $spin^c$ bundle on $Y$ with determinant line bundle $L.$ 

There is a residual gauge action of $\zz = H^1(Y; i\zz) \subset 
\Omega^1(Y; i\rr)$ on $V$ with respect to which $CSD$ is periodic modulo 
$4\pi^2 \ell;$
more precisely, the action of $u \in H^1(Y; i\zz)$ changes $CSD$ by
$4\pi^2 ([u] \cup c_1(\spc))[Y].$ The Seiberg-Witten moduli space is
therefore periodic modulo $\zz.$ It would be compact if we
divided out by the residual gauge action, but this would destroy the 
linear structure of the configuration space and hinder the application of  
finite dimensional approximation.

In \cite{M}, $\swf$ was defined as the Conley index of a flow on a finite
dimensional space. To be able to define the Conley index in the current
setting we need to restrict our attention to a bounded subset of the
configuration space. We manage to do this by cutting the configuration
space between two levels of the Chern-Simons-Dirac functional. (This
procedure was used before by Fintushel and Stern in \cite{FS} to define
$\zz$-graded instanton homology and by Marcolli and Wang in \cite{MWa} in
the context of Seiberg-Witten theory.) We approximate the gradient flow of
$CSD$ in a bounded set between the two levels by a flow on a
finite-dimensional approximation. The pro-spectrum $\swf$ is then obtained 
from the Conley index of this flow by taking direct and inverse limits as 
the levels of $CSD$ go to $\pm \infty.$

If we perturb the Seiberg-Witten equations by a nonexact $\nu$ with
$[\nu ]= - c_1(\spc),$ the Chern-Simons-Dirac functional is 
invariant under the residual gauge action, so we can no longer use it to 
cut the configuration space. Instead, we cut between the levels of $CSD + 
f,$ for a certain function $f.$ Depending on our choice of $f,$ we are 
free to take either direct or inverse limits as the levels go to $\pm 
\infty.$ Direct limits at $+ \infty$ and inverse limits at $- \infty$ give 
us back $\swf,$ while direct limits in both directions produce the new 
invariant $\swf_0.$ This is different from $\swf$: for example, when $Y= 
S^1 \times S^2$ with any nontorsion $spin^c$ structure, $\swf = *$ is 
trivial, while $\swf_0$ is nontrivial. 

We choose to mention $\swf_0$ here because of the connection to 
$\os$ theory. In \cite{OS1} and \cite{OS2}, Ozsv\'ath and Szab\'o have 
constructed several versions of Floer homology for three-manifolds, which 
they denoted by $\widehat{HF}, HF^{+}, HF^-, HF^{\infty}.$ Their theory is 
supposed
to give the same output as Seiberg-Witten theory. In \cite{OS2}, they
have made the precise conjecture relating $HF^+$ and $HF^-$ to
two versions of the Seiberg-Witten Floer homology for rational homology 
3-spheres.

More generally, we conjecture that all variants of the $\os$ Floer 
homology are different functors applied to our invariant. Given an 
integer-graded generalized homology theory $h_*$ for $S^1$-equivariant 
pro-spectra, we can apply $h_*$ to either $\swf$ or $\swf_0$ and obtain a 
sequence of abelian groups periodic modulo $\ell.$ Since a change in the 
Riemannian metric only changes $\swf, \swf_0$ by suspending or 
desuspending by an even dimensional representation,  
$h_*(\swf)$ and $h_*(\swf_0)$ are invariants of $Y$ and $\spc$ only, 
thought of either as having a relative grading modulo $\ell$ or an absolute 
grading modulo $2.$

\begin {remark} 
All generalized homology theories $h_*$ are 
meant to be reduced, but we suppress the conventional tilde from notation 
for simplicity.
\end {remark}

In section~\ref{sec:hom}, we explore some of these theories: 
the ordinary (non-equivariant) homology $H_*,$ the equivariant Borel 
homology $H^{S^1}_*,$ the co-Borel homology $cH^{S^1}_*,$ and the Tate 
homology $tH^{S^1}_*.$

\begin {conjecture}
\label {osz}
Let $Y$ be a 3-manifold with $b_1(Y)=1$ and a nontorsion $spin^c$ 
structure $\spc.$ Then:
$$  \widehat{HF}_n = H_n(\swf_0) = H_n(\swf); \hskip6pt HF^+_n = 
H^{S^1}_{n+1}(\swf_0) = H^{S^1}_{n+1}(\swf);$$
$$ HF^{-}_n = cH^{S^1}_{n+2}(\swf_0); \hskip6pt HF^{\infty}_n = 
tH^{S^1}_n(\swf_0).$$
The same should be true for manifolds with $b_1(Y)=0$ if we 
replace $\swf_0$ by $\swf.$
\end {conjecture}

In the last section we discuss relative invariants of four-manifolds with 
boundary and prove:

\begin {theorem}
\label {relative} 
Consider a Riemannian 4-manifold $X$ together with
a $spin^c$ structure, such that the boundary $(Y, \spc)$ has $b_1(Y)=1, 
c_1(\spc)$ nontorsion, and the image of $H^1(X;\rr)$ in $H^1(Y;\rr)$ is 
zero. Let 
$T(Ind)$ be the Thom spectrum for the Dirac index bundle on
$X.$ The Seiberg-Witten equations on $X$ induce
a morphism $$\Psi : \Sigma^{-b^+_2(X)} T(Ind) \to \swf(Y, 
\spc,g, A_0).$$ 
\end {theorem}

The reader may wonder what happens when $b_1(Y) \geq 1$ or $b_1(Y)=0$ and
$c_1(\spc)$ is torsion. As Cohen, Jones and Segal pointed out in
\cite{CJS}, in general there is an obstruction to defining Floer homotopy
which lies in $KO^1$ of the configuration space, and usually factors 
through the map $K^1 \to KO^1.$ For the Seiberg-Witten case, we present a 
computation of this obstruction in the appendix. It turns out that the 
obstruction is zero if and only if for every 
$$a_1, a_2, a_3 \in c_1(\spc)^{\perp} = \{a \in H^1(Y; \rr): a_j \cup 
c_1(\spc)=0\},$$ 
we have 
$$a_1 \cup a_2 \cup a_3 = 0.$$
 
For example, the obstruction is zero whenever $b_1(Y) \leq 3.$ In such 
cases, we expect to be able to mod out by the residual gauge action in the 
cohomological directions in which $CSD$ is unchanged. We can then use 
the notion of Conley index over a base from \cite{MRS}, where the base is 
a Picard torus. The result should be a ``fiberwise deforming'' 
pro-spectrum over the torus. When we vary the metric, this can change 
by suspending or desuspending by an arbitary complex bundle over the 
torus. We hope to explain this case in a subsequent paper.
  
However, the obstruction can be nonzero. This is the case, for example,
when $Y=T^3$ and $\spc$ is torsion. In such a situation, it seems that to
define a similar invariant one needs an even more general notion. Furuta
has proposed in \cite{Fu2} the concept of \emph {pro-spectrum with
parametrized universe.} It should be noted that, according to \cite{CJS},
one cannot expect to define stable Floer homotopy groups in this setting,
but at most complex oriented generalized Floer homology theories, such as
Floer complex bordism and Floer K-homology.

\section {Seiberg-Witten trajectories}
\label {sec:traj}

Let $Y$ be a closed, oriented, Riemannian 3-manifold with $b_1(Y) =1.$
Fix a nontorsion $spin^c$ structure $\spc$ on $Y$ with spinor bundle $W_0$ 
and set
$L = \det(W_0).$ We identify the space of $spin^c$ connections on $W_0$
with $i\Omega^1(Y)$ via the correspondence $A \to A - A_0,$ where $A_0$
is a fixed reference connection. We denote Clifford multiplication by
$\rho : TY \to \text { End}(W_0)$ and the Dirac operator corresponding to
the connection $A_0 + a$ by $\dir_a = \rho(a) + \dir.$

The gauge group $\gau = \text{ Map}(Y, S^1)$ acts on the space 
$i\Omega^1(Y) \oplus \Gamma(W_0)$ by
$$ u(a, \phi) = (a - u^{-1}du, u\phi).$$
We will work with the completions of $\spa$ and $\gau$ in the
$L^2_{k+1}$ and $L^2_{k+2}$ norms, respectively, where $k\geq 4$ is a fixed
integer. In general, we denote the $L^2_k$ completion of a space $E$ by
$L^2_k(E).$

Unlike in \cite{M}, here it becomes necessary to perturb the Seiberg-Witten 
equations by 
an exact 2-form $\nu$ on $Y$ in order to obtain a genericity 
condition. The perturbed 
Chern-Simons-Dirac functional is defined on $L^2_{k+1}(\spa)$ as
$$ CSD_{\nu}(a, \phi) = -\frac{1}{2} \int_Y a \wedge (2F_{A_0}+ 
da+2\pi i\nu) 
+ \int_Y \langle \phi, \dir_a \phi \rangle \dif \text{vol} .$$
The change in $CSD_{\nu}$ under the action of the gauge group is
$$CSD_{\nu} (u\cdot (a, \phi)) - CSD_{\nu}(a, \phi) = 4\pi^2 ([u] \cup 
c_1(\spc))[Y].$$

The gradient of $CSD_{\nu}$ with respect to the $L^2$ metric 
is the vector field
$$ \nabla CSD_{\nu} (a, \phi) = (*da +  *F_{A_0} + i\nu +\tau(\phi, 
\phi), \dir_a \phi),$$
where $\tau$ is the bilinear form defined by $\tau(\phi, \psi) = 
\rho^{-1}(\phi \psi^*)_0$ and the subscript $0$ denotes the trace-free part.

The \emph {perturbed Seiberg-Witten equations} on $Y$ are the equations for 
the critical points of $CSD_{\nu}:$
$$ *da + *F_{A_0} + i\nu +\tau(\phi, \phi)=0, \hskip 6pt \dir_a \phi 
=0.$$

The basic compactness result for the solutions $(a, \phi)$ to 
the Seiberg-Witten equations (\cite{K}) is that one can always find a gauge 
transformation $u$ such that $u(a, \phi)$ is smooth and bounded in all 
$C^m$ norms by constants which depend only on $\nu$ and the metric on 
$Y.$

As mentioned in the introduction, we intend to cut the moduli space
between two levels of $CSD_{\nu}.$ In order
for this to be possible, we need to impose a 
constraint on $\nu.$ 

\begin {definition}
\label {def}
A perturbation $\nu$ is called {\bf good} if the set of
critical points of $CSD_{\nu}$ is discrete modulo gauge. 
\end {definition}

The set of good
perturbations is Baire in the space of exact 2-forms. Indeed, according to
\cite{Fr}, all critical points are nondegenerate for a Baire set of
perturbations, and any nondegenerate critical point is isolated.

From now on we will always assume that $\nu$ is good. Let us study the 
trajectories 
of the downward gradient flow of $CSD_{\nu},$ given by:
\begin {equation}
\label{csdi}
x = (a, \phi): \rr \to L^2_{k+1}(\spa), \delt x(t) = -\nabla CSD_{\nu} 
(x(t)).   
\end {equation}

Such Seiberg-Witten trajectories can be interpreted as solutions of the 
four-dimensional Seiberg-Witten equations on the infinite cylinder $\rr \times 
Y.$

We borrow the terminology of \cite{M} and say that a Seiberg-Witten trajectory 
$x(t)$ is \emph {of finite type} if
both $CSD_{\nu}(x(t))$ and $\|\phi(t)\|_{C^0}$ are bounded functions of 
$t.$

The following proposition was proved in \cite{M} for the case $b_1(Y)=0,$ 
but the proof works in general, with only minor changes:

\begin {proposition}
\label {trajectories}
There exist $C_m > 0$ such that for any $(a,\phi)\in L^2_{k+1}(\spa)$
which is equal to $x(t_0)$ for some
$t_0\in \rr$ and some Seiberg-Witten trajectory of finite type $x:\rr\to
L^2_{k+1}(\spa)$, there exists $(a', \phi')$ smooth and gauge equivalent
to $(a, \phi)$ such that $\|(a', \phi')\|_{C^{m}} \leq C_m$ for all $m > 0.$
\end {proposition}

As in \cite{M}, it is useful to project our configuration space to the 
\emph{Coulomb gauge slice}. Let $\go$ be the group of gauge 
transformations of the form $u=e^{i\xi},$ where $\xi:Y \to \rr$ 
satisfies $\int_Y \xi =0.$ The Coulomb gauge slice is the space
$$V = i\ker d^* \oplus \Gamma(W_0) \cong \spa/\go.$$

For every $(a, \phi) \in \spa,$ there is a unique element $\Pi(a, \phi)\in 
V$ which is gauge equivalent to $(a, \phi)$ by a transformation in $\go.$ 
We call it the \emph{Coulomb projection} of $(a, \phi).$ 

There is a residual gauge action of $H^1(Y; \zz) \times S^1$ on $V$ as 
follows: if we choose basepoints on $Y$ and $S^1,$ then $h \in H^1(Y; 
\zz)$ is the homotopy class of a unique pointed harmonic map $u : Y \to 
S^1.$ Then $h$ acts on $(a, \phi)$ via the gauge transformation $u,$ 
while $e^{i\theta} \in S^1$ sends $(a, \phi) \in V$ to $(a, 
e^{i\theta}\phi).$  

As in \cite{M}, one can find a metric on $V$ such that the downward 
gradient trajectories of $CSD_{\nu}|_V$ with respect to this metric are 
exactly the Coulomb projections of the trajectories of $CSD_{\nu}$ on 
$\spa.$ Given a tangent vector at some point in $V,$ its length in the 
new metric is the $L^2$ length of its projection to the orthogonal 
complement of the tangent space to the gauge equivalence class through 
that point in $\spa.$ The gradient of $CSD_{\nu}|_V$ with respect to this 
metric on $V$ is of the form $l + c,$ where $l, c: V \to V$ are given by
$$ l(a, \phi) =(*da, \dir \phi)$$
$$ c(a, \phi) = (\pi \circ (F_{A_0} +\tau(\phi, \phi)) + i\nu, 
\rho(a)\phi - i\xi(\phi)\phi).$$

Here $\pi$ denotes the orthogonal projection from $\Omega^1(Y)$ to $\ker
d^*.$ Also, $\xi(\phi):Y \to \rr$ is defined by $d\xi(\phi) =
i(1-\pi)\circ \tau(\phi, \phi)$ and $\int_Y \xi(\phi) =0.$ 

Let us look at finite type trajectories $x: \rr \to L^2_{k+1}(V)$ for 
some fixed $k \geq 4.$ From Proposition \ref{trajectories} we know that 
they are locally the Coulomb projections of smooth trajectories contained 
in a bounded set modulo the residual gauge action. In other words, if we 
denote by $Str_m(R)$ the union of balls
$$ \{ (a, \phi) \in V | \exists h \in H^1(Y; \zz) \text{ such that }
\|h \cdot (a, \phi) \|_{L^2_m} \leq R\},$$
then the following statement is true: all trajectories $x$ as above are 
smooth in $t$ and there are uniform constants $C_m > 0$ such that 
$x(t) \in Str_m(C_m)$ for each $m.$

\smallskip

\section {The construction of $\swf$}
\label {sec:g}

This section contains the proof of Theorem~\ref{nontorsion} (a). Let $h$
be the generator of $H^1(Y; \zz)$ which satisfies $$ \ell = (h \cup
c_1(\spc))[Y] > 0.$$

Let $u: Y\to S^1$ be the pointed harmonic map in the homotopy class $h.$ 
Then $u^{-1}du = ih,$ where we think of $h$ as a harmonic 1-form. 
The Seiberg-Witten moduli space is compact modulo the residual 
gauge action $u \cdot (a, \phi) \to (a - ih, u\phi).$ We have 
$CSD_{\nu}(u\cdot (a, \phi)) = CSD_{\nu} (a, \phi) + 4\pi^2 \ell.$

We follow \cite{M} and consider the orthogonal projections from $V$ to 
the finite dimensional subspaces $V^{\mu}_{\lambda}$ spanned by the 
eigenvectors of $l$ with eigenvalues in the interval $(\lambda, \mu].$ We 
can smooth out these projections to obtain a family $p^{\mu}_{\lambda}$ 
which is continuous in $\mu$ and $\lambda$ and still satisfies 
$V^{\mu}_{\lambda} = \im (p^{\mu}_{\lambda}).$ 

\subsection {The Conley index.} As mentioned in the introduction, we 
intend to cut the moduli space 
between two levels of the $CSD_{\nu}$ functional. Since the set of 
critical points is discrete modulo gauge, we can choose $v \in \rr$ which 
is not the 
value of $CSD_{\nu}$ at any critical point. Because of the periodicity in 
the residual gauge direction, the same must be true for the values $v + 
4\pi^2 n\ell, n \in \zz.$ Choose $m, n \in \zz, m<n,$ and consider the 
set 
$$ T(R) = \{(a,\phi) \in Str_{k+1}(R) | CSD_{\nu}(a, \phi) \in (m,n)\}.$$

Recall that the ``strip'' $Str_{k+1}(R)$ is the union of the residual
gauge translates of the ball of radius $R$ in the $L^2_{k+1}$ norm. On
that ball, $CSD_{\nu}$ takes values in a compact interval $I \subset \rr.$
On a translate of the ball, it takes values in the corresponding interval
$I + 4\pi^2n\ell,$ for some $n \in \zz.$ It follows that $T(R)$ intersects
only finitely many such translates, hence it is bounded in the $L^2_{k+1}$
norm.  Furthermore, if $R$ is sufficiently large, all the Seiberg-Witten
trajectories contained in $T(2R)$ are of finite type, hence contained in
$T(R).$

Now we are in the right setting for doing finite dimensional 
approximation: the gradient flow of $CSD_{\nu}$ on the bounded set 
$T(2R)$ is of the form $-\delp x(t) = (l+c)x(t),$ where $l$ is linear 
Fredholm and self-adjoint and $c: L^2_{k+1}(V) \to L^2_k(V)$ is compact.
Let us consider the trajectories of the gradient flow of $CSD_{\nu}$ 
restricted to $V^{\mu}_{\lambda}$ which are contained in $T(2R).$ The 
following compactness result is the analogue of Proposition 2 in 
\cite{M}, and the proof is completely similar:

\begin {proposition}
\label {main}
For any $-\lambda$ and $\mu$ sufficiently large, if a trajectory $x: \rr
\to L^2_{k+1}(V_{\lambda}^{\mu}),$
$$ (l+p_{\lambda}^{\mu} c)(x(t)) = -\delt x(t) $$
satisfies $x(t) \in T(2R)$ for all $t,$ then in fact $x(t) \in
T(R)$ for all $t.$
\end {proposition}

Let $S$ be the invariant part of $T=T(2R) \cap V^{\mu}_{\lambda}$ under
the flow, i.e. the set of critical points of 
$CSD_{\nu}|_{V^{\mu}_{\lambda}}$ 
contained in $T$ together with the gradient trajectories between them. 
Then $S$ is contained in the interior of $T$ by Proposition \ref{main} 
and the fact that no gradient trajectory can be tangent to a level set  
of $CSD_{\nu}.$

Because of these properties, we can associate to $S$ a \emph{Conley 
index} $I(S),$ which is the pointed space $N/L,$ where $(N, L)$ is an 
index pair for $S,$ i.e. a pair of compact subsets $L \subset N \subset 
T$ satisfying the following conditions:

\begin {enumerate}
\item $S \subset \inte (N \setminus L);$
\item $L$ is an exit set for $N,$ i.e. for any $x\in N$ and $t>0$ such
that $\vp_t(x) \not\in N,$ there exists $\tau \in [0,t)$ with
$\vp_{\tau}(x) \in L.$ Here we denote by $\phi$ the downward gradient 
flow on $T.$

\item $L$ is positively invariant in $N,$ i.e. if for $x\in L$ and $t>0$
we have $\vp_{[0,t]}(x) \subset N,$ then in fact $\vp_{[0,t]}(x) \subset
L.$
\end {enumerate}

For the basics of Conley index theory, the reader is referred to 
\cite{C} or \cite{Sa}. Section 5 of \cite{M} also gives an overview of 
the relevant properties. The most important ones are the existence of the 
index pair and the fact that the Conley index is independent of $N$ and 
$L$ up to canonical homotopy equivalence.

We are interested in the Conley index $I^{\mu}_{\lambda}(m,n)$ of $S = 
S^{\mu}_{\lambda}(m,n)$ because its homology is the 
same as the Morse homology computed by counting critical points and 
trajectories between them in the usual way. In our case, if we are able 
to take the limits $n, \mu \to \infty, m, \lambda \to -\infty,$ the 
homology of the resulting object should be some version of Seiberg-Witten 
Floer homology. 

Let us first see what happens as $\mu \to \infty$ and $\lambda \to
-\infty.$ This process was studied in detail in \cite{M}. If $\mu$
increases so that $V^{\mu}_{\lambda}$ increases in dimension by 1, the
flow on the new $T$ is homotopic to the product of the flow on the old
$T$ and a linear flow on the complementary subspace. Since the linear
operator $l$ on this subspace has only positive eigenvalues, the
respective Conley index is trivial, which implies that
$I^{\mu}_{\lambda}$ does not change with $\mu,$ up to homotopy
equivalence.

However, if we decrease $\lambda$ so that $V^{\mu}_{\lambda}$ increases
in dimension by 1, then $l$ has negative eigenvalues on the complementary
subspaces and the new Conley index $I^{\mu}_{\lambda}$ is the suspension
of the old one. In order to obtain an invariant object, we need to
desuspend by $V^0_{\lambda}.$ Let $\mathcal{S}(I^{\mu}_{\lambda})$ be the $S^1$-equivariant 
spectrum of $I^{\mu}_{\lambda}$, in the sense of \cite{LMS}.
Set:
$$ J^{\mu}_{\lambda}(m,n) = \Sigma^{-V^0_{\lambda}} \mathcal{S}(I^{\mu}_{\lambda}).$$ 
Then $J(m,n) =J^{\mu}_{\lambda}(m,n)$ is independent of $\mu$ and 
$\lambda,$ up to canonical equivalence.

\subsection {Some algebraic-topological preliminaries.}

Let $\sp$ be the homotopy category of $S^1$-equivariant spectra 
with semifree $S^1$ actions, modelled 
on the standard universe $\rr^{\infty} \oplus \cc^{\infty}.$ Here $\rr$ 
and $\cc$ are the real and complex representations of $S^1.$ 

As it stands, the category $\sp$ does not have colimits. However, as 
explained in \cite{Ma} for the nonequivariant case, given a sequential 
direct system of spectra in $\sp:$
$$ X_1 \to X_2 \to \cdots $$
there is a \emph {minimal weak colimit} $X=$ wcolim $X_i$ with maps $X_i 
\to X$ inducing isomorphisms at the level of all homotopy groups. The 
minimal weak colimit is not usually a colimit: given a system of commuting 
morphisms $X_i \to Y,$ they factor through a morphism $X \to Y,$ but this 
morphism is not unique. Furthermore, while $X=$ wcolim $X_i$ is unique up 
to equivalence, it is not unique up to canonical equivalence. Indeed, 
there can be nontrivial self-homotopy equivalences of $X$ which are the 
identity on all homotopy groups: they could be of the form id $+f,$ where 
$f$ is a \emph{phantom map}.

\begin {definition}
A {\bf phantom map} $f: X \to Y$ is a morphism in $\sp$ with the property 
that for every finite CW-spectrum $Z \in \text{Ob }\sp$ and morphism $g: Z 
\to X,$ the composite $f \circ g: Z \to Y$ is zero. We denote the
abelian group of phantom maps from $X$ to $Y$ by ph$(X,Y).$
\end {definition}

Note that this maps are called {\bf f}-phantoms in \cite{Ma}. 
We can get a better behavior if we replace the category $\sp$ by $\sp',$ 
whose objects are the same as those of $\sp,$ but 
whose sets of morphisms 
are:
$$ \text{Mor}_{\sp'}(X, Y) = \text{Mor}_{\sp}(X, Y)/ \text{ph}(X, Y).$$

From now on we will denote $ \text{Mor}_{\sp'}(X, Y)$ by $[X,Y]'_{S^1}$ 
for simplicity. The weak colimit becomes an actual colimit in the category 
$\sp'.$ In other words, for a sequential direct system $X_i$ as before, 
there exists $X,$ unique up to canonical equivalence, such that for every 
$Y \in \sp$:
$$ [X, Y]'_{S^1} = \lim_{\from} [X_i, Y]'_{S^1}.$$  

Moreover, for every generalized homology theory $h_*,$ there is a natural 
isomorphism:
$$ h_*(X) = \lim_{\to} h_*(X_i).$$

The category $\sp'$ still does not have good inverse limits. There is a 
notion of weak inverse limit, but its cohomology is not the direct limit 
of the cohomologies of the terms. 

The solution to this problem is to introduce pro-spectra, along the lines
of Cohen, Jones, and Segal (\cite{CJS}). (For the general definition of a
pro-category, we refer to \cite{AM}.) Basically, a \emph {pro-spectrum}
$X$ is an inverse system of spectra $\{ X_p\}, p\in \zz$:
$$ \cdots \from X_{p-1} \from X_p \from \cdots $$

We could call $X$ the \emph {pro-limit} of $X_p.$ The set of morphisms 
between two pro-spectra $X = \{X_p\}$ and $Y = \{Y_q\}$ is defined as
$$ \text{Mor}_{\pros}(X, Y) = \lim_{\infty \from q} \lim_{p \to \infty}  
[X_p, Y_q].$$

In fact, it is easy to check that the pro-limit is in fact the inverse limit
of $X_p$ in the category $\pros$ of pro-spectra. As before, we 
denote $\text{Mor}_{\pros}(X, Y)$ by $[X, Y]'_{S^1}$ for simplicity.
The category $\pros$ is closed under both sequential colimits and limits.

\subsection {Taking co- and pro-limits.} 
Let us come back to the spectrum $J(m,n),$ which still depends on $m$ and
$n,$ and let us see what happens as we vary $n.$ In a finite 
dimensional
approximation $V^{\mu}_{\lambda},$ the isolated invariant set
$S^{\mu}_{\lambda}(m,n)$ is an attractor subset of $S=S^{\mu}_{\lambda} 
(m,n+1).$ This means that it ``attracts'' nearby points
of $S$ under the downward gradient flow, which is obviously true because
$CSD_{\nu}$ is decreasing along flow lines. There is a corresponding
repeller subset $S^{\mu}_{\lambda}(n,n+1).$ There is a coexact sequence
for the Conley indices of an attractor-repeller pair (see \cite{C}):
$$ I(S^{\mu}_{\lambda}(m,n)) \to I(S^{\mu}_{\lambda}(m,n+1)) \to
I(S^{\mu}_{\lambda}(n,n+1)) \to \Sigma I(S^{\mu}_{\lambda}(m,n)) \to 
\dots$$

For $\mu$ and $-\lambda$ sufficiently large, the first map in the above 
sequence gives a morphism $J(m,n) \to J(m, n+1).$ This gives a sequential 
direct system:
$$ J(m, n) \to J(m, n+1) \to J(m, n+2) \to \dots.$$
we can take its colimit in $\sp'$:
$$ J(m, \infty) = \colim_n J(m,n).$$ 
 
We may expect a similar construction for $m \to -\infty.$ Indeed, we
have exact triangles: 
$$ J(m-1,m) \to J(m-1, n) \to J(m,n) \to \Sigma J(m-1,m) \to \dots $$ 

The maps obtained by the composition $J(m-1, n) \to J(m, n) \to J(m, 
\infty)$ induce a well-defined morphism $J(m-1, \infty) \to J(m, \infty)$ in 
$\sp.$ 

Using the resulting maps we can take the pro-limit of $J(m, 
\infty)$ as $m \to -\infty.$ Let us define $SWF(Y, \spc, g)$ as being the 
pro-spectrum 
$J(-\infty, \infty) = \{J(m, \infty)\}.$

\subsection {Invariance.}
Of course, we have to check that our invariant is independent of the
choices made in its construction. The first step is:

\begin {proposition}
\label {inv}
The pro-spectrum $J(-\infty, \infty)$ does not depend on the
value $v$ of $CSD_{\nu}$ where we do the cutting, up to canonical 
equivalence in $\pros.$ 
\end {proposition}

\noindent{\bf Proof. } From our construction it
is clear that nothing changes if we replace $v$ by $v+4\pi^2n\ell, n\in
\zz$ (the direct limits and the pro-limits are the same). Thus it suffices
to study the case of $v, v' \in \rr$ which are not values of $CSD_{\nu}$
at critical points and which satisfy $v < v' < v + 4\pi^2\ell.$ Let us
switch notation and denote by $J(a,b)$ the Conley index obtained as before
from the approximate trajectories between the levels $CSD_{\nu} =a$ and
$CSD_{\nu} =b,$ for some $\mu, -\lambda \gg 0.$ (For example,
$J(v+4\pi^2m\ell, v+4\pi^2n\ell)$ is what we previously denoted $J(m,n).$)
We have a sequential direct system:
$$ J(v+4\pi^2m\ell, v+4\pi^2n\ell) \to J(v+4\pi^2m\ell, v'+4\pi^2n\ell) 
\to $$ 
$$ \to J(v+4\pi^2m\ell, v+4\pi^2(n+1)\ell) \to J(v+4\pi^2m\ell, 
v'+4\pi^2(n+1)\ell) \to \dots$$

Its colimit is the same as that of its subsequence 
$$ J(v+4\pi^2m\ell, v+4\pi^2n\ell) \to J(v+4\pi^2m\ell, 
v+4\pi^2(n+1)\ell)
\to \dots$$ which gives back
$J(v+4\pi^2m\ell, \infty),$ the one previously denoted $J(m, \infty).$

Taking the pro-limits as $m \to -\infty$ according to the levels of both
$v + 4\pi^2 \zz$ and $v'+4\pi^2 \zz$ gives the same pro-spectrum as taking
the pro-limit according to each of them separately (they would be 
subsystems of the same inverse system.) $\fin$

\medskip
Then, from the invariance of the Conley index under deformations 
it is straightforward to deduce:

\begin {proposition}
$\swf(Y, \spc, g, A_0)$ does not change (up
to canonical equivalence) as we vary the other parameters involved in the
construction, such as the perturbation $\nu$ or the radius $R$ in
Proposition \ref{main}. 
\end {proposition}

\subsection{Changing the connection and the metric.} Let us explain what 
happens when we vary the base connection $A_0.$ It suffices to consider 
nearby connections $A_0, A_0'.$ For these we can find $\mu$ and $-\lambda$ 
sufficiently large so that $V^{\mu}_{\lambda}$ does not changes dimension 
as we choose the base connection to be $A_0 + t(A_0' - A_0), t\in [0,1].$ 
Using the invariance under continuation of the Conley index, we get that 
all the corresponding Conley indices are equivalent. However, the number 
of negative eigenvalues $n_{\lambda}=\dim V^0_{\lambda}$ by which we 
desuspend in the construction of $\swf$ varies according to the spectral 
flow of the Dirac operator $\dir_t, t\in [0,1].$ Let $E$ be the spinorial 
part of $V^0_{\lambda}$ for $\dir_0$ and $E'$ that for $\dir_0.$ Then $E, 
E'$ are complex $S^1$-representations, and the difference in their 
dimensions is the spectral flow. After taking the 
co- and pro-limits, it follows that $\swf$ changes according to the 
formula:
$$ \swf(Y, \spc, g, A_0') = \Sigma^{E' - E} \swf(Y, \spc, g, A_0).$$

Changes in the Riemannian metric have a similar effect.  

\subsection {Periodicity.}
Let us consider a particular change in connection: the homotopy $A_t = A_0 
- ith, t\in [0,1].$ The spectral flow of the Dirac operators (over 
$\rr$) along this homotopy is $\ell.$ Hence:
\begin {equation}
\label{unu}
(I^{\mu}_{\lambda} (m,n) \text{ with } A_0) \cong \Sigma^{E} 
(I^{\mu}_{\lambda} (m,n) \text{ with } A_0 - ih), 
\end {equation}
where $E$ is (noncanonically) isomorphic to $\cc^{\ell/2}.$

On the other hand, recall the periodicity on $V: CSD_{\nu}(a-ih, 
u\phi) = CSD_{\nu} (a, \phi) + 4\pi^2 \ell.$ This does not translate into 
a periodicity on $V^{\mu}_{\lambda}$ because $(a-ih, u\phi)$ may not be in 
$V^{\mu}_{\lambda}$ for $(a, \phi) \in V^{\mu, \lambda}.$ However, if $(a, 
\phi) \in V^{\mu}_{\lambda}$ with base connection $A_0,$ then it is true 
that $(a-ih, u\phi) \in V^{\mu}_{\lambda}$ with base connection $A_0 - 
ih.$ Hence there is a canonical equivalence:
\begin {equation}
\label{doi}
(I^{\mu}_{\lambda} (m,n) \text{ with } A_0) \cong (I^{\mu}_{\lambda} 
(m+1,n+1) \text{ with } A_0 - ih). 
\end {equation}

Putting (\ref{unu}) and (\ref {doi}) together and taking the limits as $n,
-m \to \infty$ we get the equivalence: $$ \Sigma^{E}(\swf(Y, \spc))
\cong \swf(Y, \spc)$$
mentioned in Theorem~\ref{nontorsion} (a).
\medskip

\subsection{The $S^1$ action.} All the constructions above can be
done in an $S^1$-equivariant manner, preserving the residual gauge action
of $S^1$ on $V,$ given by $e^{i\theta}: (a, \phi)\to (a, 
e^{i\theta}\phi).$
The resulting pro-spectrum $\swf$ is $S^1$-equivariant.

In fact, in this case we can ``divide
out'' by the $S^1$ action. Indeed, there is no reducible critical point 
$(a,0)$ of $CSD_{\nu},$ because there are no flat connections on $W_0.$
This is also true for critical points in the finite dimensional
approximations. It follows by $S^1$-equivariance that the flow lines
between such points also do not intersect the plane $\phi =0.$ Thus we
can replace the strip $Str$ by the set $Str'$ obtained from $Str$ by  
deleting a neighborhood of the plane $\phi=0.$ The $S^1$ action is free
on $Str',$ so we can take the quotient $Str''=Str'/S^1.$ By doing all the
constructions as before, but with the quotient flow on $Str'',$ we obtain
a pro-spectrum $\swfn.$ The periodicity is reflected in an equivalence 
$\Sigma^{\ell} \swfn \cong \swfn.$
\medskip

\subsection {Duality}

Every $S^1$-equivariant spectrum has a Spanier-Whitehead dual defined as 
the function spectrum:
$$ DX = F(X, \ss),$$
where $\ss = S^0$ is the sphere spectrum. It is characterized by the property 
that 
$$ [W, DX]_{S^1} = [W \wedge X, \ss]_{S^1}$$ 
for every spectrum $W.$ Furthermore, the same is true in $\sp':$
$$ [W, DX]'_{S^1} = [W \wedge X, \ss]'_{S^1}.$$

The functor $D$ is contravariant. We can extend its definition to 
pro-spectra $X =\{ X_1 \from X_2 \from \dots \}$ by letting:
$$ DX = \colim \hskip3pt DX_i.$$

If we have a sequential direct system of pro-spectra $X_1 \to X_2 \to 
\cdots,$ it is not hard to check that
$$ D(\colim X_i) = \text{pro-lim }DX_i.$$
\smallskip

\noindent {\bf Proof of Theorem~\ref{nontorsion} (c).} If we change the 
orientation on $Y,$ the function $CSD_{\nu}$ switches sign and the 
Seiberg-Witten flow changes its direction. As a consequence, by a duality 
theorem for Conley indices (see \cite{Co}), the spaces 
$J(m,n)$ for $X$ and $J(-n, -m)$ for $\bar X$ are duals. Since pro-limits 
and colimits are dual operations, it follows that $\swf= J(-\infty, 
\infty)$ for $X$ and $\swf$ for $\bar X$ are duals.

The fact that the analogue is not true for $\swf_0$ can be seen from 
the example of $S^1 \times S^2$ in Section~\ref{sec:hom}.\hskip2pt $\fin$
\medskip

\section {Nonexact perturbations}

\subsection {Reconstructiong $\swf$} In Section~\ref{sec:traj} we have 
considered exact perturbations $\nu$ of 
the $CSD$ functional. More generally, let us consider 
$$  CSD_{\nu}(a, \phi) = -\frac{1}{2} \int_Y a \wedge (2F_{A_0}+
da+2\pi i\nu)
+ \int_Y \langle \phi, \dir_a \phi \rangle \dif \text{vol} $$   
for any 1-form $\nu.$ The change under the gauge group is:
$$CSD_{\nu} (u\cdot (a, \phi)) - CSD_{\nu}(a, \phi) = 4\pi^2 ([u] \cup
(c_1(\spc) + [\nu])[Y].$$

An important qualitative difference appears when $[\nu] = -c_1(\spc):$ no
perturbation $\nu$ is good (in the sense of Definition~\ref{def}) in this
cohomology class! Indeed, there is always a line of reducible monopoles.
Thus we need to introduce an additional perturbation. Recall that $h$ is  
the generator of $H^1(Y; \zz)$ which satisfies $ \ell = (h \cup
c_1(\spc))[Y] > 0.$ We can think of $h$ as a harmonic function.
 
We replace $CSD_{\nu}$ by the functional
$$ CSD_{\nu, \epsilon}(a, \phi) = CSD_{\nu}(a, \phi) + \epsilon 
\sin (2\pi \langle a, h\rangle). $$
The (twice perturbed) Seiberg-Witten equations $\nabla CSD_{\nu, \epsilon} 
(a,\phi) =0$ have a discrete set of solutions for generic $\nu$ and small
$\epsilon.$ We call the pair $(\nu, \epsilon)$ \emph {good} if this 
condition is satisfied.

However, even introducing $\epsilon$ does not make the situation 
completely analogous to that in the previous section. The problem is that 
$CSD_{\nu}$ is periodic:
$$CSD_{\nu}(a-ih, u\phi) = CSD_{\nu}(a, \phi).$$
This implies that cutting the configuration space between its 
levels does not work as well: the set $T(R) = Str_{k+1}(R) \cap \{(a,
\phi) \in V| CSD_{\nu}(a, \phi) \in (m,n)\}$ may not be bounded.
Therefore, we need to find another way of cutting. For this we introduce 
a new notion, that of transverse functions.

The Hodge decomposition of 1-forms gives $V=ih\rr \oplus i\im d^* \oplus 
\Gamma(W_0).$ Let $p: V \to ih\rr \cong \rr$ be the orthogonal projection 
to the first factor, and $T:V\to V, T(a, \phi) = (a+ih, u^{-1}\phi)$ be 
the generator of the residual gauge group. The strip of balls $Str = 
Str_{k+1}(2R)$ is as in the previous section.
\medskip

\begin {definition} A {\bf positively (resp. negatively) transverse 
function} is a smooth functions $f : Str \to \rr$  satisfying the following 
properties:

\begin {enumerate}
\item There exists a constant $M >0$ such that $f(x) < 0$ whenever $p(x) 
< -M$ and $f(x) > 0$ whenever $p(x) > M.$ 

\item If $f(x) \geq 0,$ then $f(Tx) \geq 0.$

\item We have $\langle \nabla CSD_{\nu, \epsilon}(x), \nabla f(x) \rangle 
>0$ (resp. $< 0$ whenever $f(x) =0.$
\end {enumerate}
\end {definition}

The inner product in condition 3 is the one used for getting the gradient
of $CSD_{\nu, \epsilon}$ in the gauge slice. Condition 3 basically says
that the level sets of $f$ at $0$ are transverse to the gradient flow, and
specifies the direction of the flow at these level sets. Note that because
$\nabla CSD_{\nu, \epsilon}$ is invariant under $T,$ all the translates 
$T^n
\{x|f(x) =0 \}$ are also transverse to the flow.

Given a positively transverse function $f,$ we obtain a nice 
partition of the strip $Str$ in the following way. Let us denote 
$$A_n = T^n\{x \in Str| f(x) \leq 0\}.$$

Because of Condition 2 in the definition of a transverse function, we 
have a nested sequence 
$$ \dots \subset A_n \subset A_{n+1} \subset 
\dots$$ 
Let $U(m,n) = A_n \setminus A_m.$ 

We claim that $U(m,n)$ is bounded in the $L^2_{k+1}$ norm. Indeed,
$Str$ is the union of the residual gauge translates of a ball, so it
suffices to check that $p(U(m,n))$ is bounded in $\rr.$ But this is 
true because of Condition 1. 

Note that when $\epsilon =0,$ the function $CSD_{\nu}$ itself is 
positively transverse. In general we have:

\begin {lemma}
\label {transverse}
For every good pair $(\nu, \epsilon)$ with $[\nu] = tc_1(\spc), t \geq 
-1,$ positively transverse functions exist. When $t=-1,$ negatively 
tranverse functions also exist.
\end {lemma}

\noindent {\bf Proof.} 
Let $[\alpha, \beta] \subset (0,1)$ be a small interval such that $\nabla 
CSD_{\nu, \epsilon} (x) \neq 
0$ for all $x \in Str$ such that $p(x) \in (\alpha, \beta).$ We then define a 
smooth increasing function $g: [0,1] \to \rr$ such that $g \equiv 0$ on 
$[0, \alpha],$ and $f$ increases slowly on $[\alpha, \beta]$ so that 
$$ \|\nabla (g\circ p)(x)\|  < \|\nabla CSD_{\nu, \epsilon} (x) \| \text { 
if } x \in Str, p(x) \in [\alpha, \beta].$$
Then we set $g\equiv \delta$ on $[\beta, 1]$ for some $\delta = 
g(\beta) > 0.$ Next, we extend $f$ to the whole real line by requiring 
$g(u+1) = g(u) + \epsilon$ for all $u \in \rr.$ Set :
$$ f = g\circ p + CSD_{\nu, \epsilon}.$$

We claim that $f$ is positively transverse. Condition 1 is satisfied
because $\lim_{u \to \pm\infty} g(u) = \pm\infty$ and $CSD_{\nu, 
\epsilon}$ is either invariant under $T$ and therefore a bounded function 
on $Str$ (when $t = -1$) or satisfies $\lim_{u \to \pm\infty} CSD_{\nu, 
\epsilon}(u) = \pm\infty$ itself (when $t \geq -1.$) 

Condition 2 is also
satisfied because 
$$f(Tx) \geq f(x) + \delta$$ for $x \in
Str$ (when $t=-1$ we have equality.) 

Finally, condition 3 is satisfied  because $\|\nabla (g\circ p)\|
\leq \|CSD_{\nu, \epsilon}\|,$ with equality only at the critical points 
of $CSD_{\nu, \epsilon},$ and we can easily arrange so that that the level 
sets of $f$ at $0$ do not go through such points. 

Similarly, for $t = -1$ one can show that the function $g\circ p - 
CSD_{\nu, \epsilon}$ is negatively transverse. However, this function does 
not satisfy conditions 1 and 2 in the case $t \geq -1,$ when $CSD_{\nu, 
\epsilon}$ is not periodic. $\fin$

\medskip

Now that we have constructed a positively transverse function, note that
each of the sets $U(m,n)$ can be used instead of $T(2R)$ to do finite
dimensional approximation as in the previous section. Indeed, they are
bounded and the flow lines are transverse to the level sets which are
separating them, because of condition 3.

We can again consider the set of critical points of $CSD_{\nu, \epsilon}$ 
inside of $U_n$ together with the trajectories between them (in a suitable finite 
dimensional approximation), and take its Conley index. After taking the 
relevant desuspension, we get spectra $J^i_n = J(U_n),$ the analogues of 
$J(m,n)$ from the previous section. We can take colimits as $n 
\to \infty$ and pro-limits as $m \to -\infty.$ The result is a 
pro-spectrum $J(-\infty, \infty).$ 

\begin {proposition}
The pro-spectrum $J(-\infty, \infty)= \swf(Y, \spc, g, A_0)$ is 
independent of $\epsilon, \nu$ 
and of the positively transverse function $f$ used in its definition, as 
long as $[\nu] = tc_1(\spc)$ for $t \geq -1,$ and $\epsilon \neq 0$ when 
$t=-1.$ 
\end {proposition}

\noindent{\bf Proof. } The proof of invariance under changing the 
transverse function $f$ is similar to the proof of Proposition~\ref{inv}.
Independence of the perturbation follows from the invariance under 
continuation of the Conley index. $\fin$

\medskip

\subsection {The construction of $\swf_0$} \label {sec:s0} For 
$[\nu]=-c_1(\spc),$ 
there is an alternate construction, which makes 
use of negatively transverse functions as well. Let $f_1$ be a positively 
transverse function and $f_2$ a negatively transverse one. We denote:
$$A_n = T^n\{x \in Str| f_1(x) \leq 0\}; \hskip6pt B_n = T^n\{x \in Str| 
f_2(x) \leq 0\}.$$

For $n \gg 0 \gg m,$ we have $B_m \subset A_n.$ Let $V(m,n) = A_n 
\setminus B_m.$ There is a nested sequence:
$$ V(m, n) \subset V(m-1, n+1) \subset V(m-2, n+2) \subset \dots$$

If we take the Conley indices of the flow inside $V(m,n)$ and desuspend by 
$V^0_{\lambda}$ we get finite spectra $J'(m,n).$ Note that $V(m,n)$ is 
an attractor subset of $V(m-1, n+1).$ Thus we can take the 
colimit as 
$n \to \infty, m\to -\infty.$ Set  
$$\swf_0(Y, \spc, g, A_0) = \colim_{m,n} J'(m,n).$$ 
 
Unlike $\swf,$ since there is no need of taking pro-limits, this is an 
actual spectrum. Apart from this, it is easy to see that the other 
invariance properties of $\swf$ still hold for $\swf_0.$ 

Let us now construct the morphism $j: \swf_0 \to \swf$ announced in the 
statement of Theorem~\ref{nontorsion}. Let $U(m,n)$ be the sets 
between the different levels of $f_1$ as before. Then, for $m \ll m' \ll 
0,$ $V(m',n)$ is an attractor subset of $U(m,n).$ This induces a map 
between the Conley indices:
$$ J'(m', n) \to J(m,n).$$
Sending in turn $n \to \infty, m \to -\infty,$ and $m' \to -\infty,$ we 
obtain the desired morphism in $\pros$:
$$ j: J'(-\infty, \infty) \to J(-\infty, \infty).$$

It is not hard to check that this does not depend on the choices made in 
its construction.

\section {Floer Homologies}
\label{sec:hom}

Let $X$ be an $S^1$-equivariant pointed pro-spectrum and $DX$ its dual. In 
this section we discuss some of the generalized homology and cohomology 
theories of $X,$ and what happens when we apply them to $\swf$ and 
$\swf_0.$ All our theories are reduced, but we do not write down the 
conventional tilde.

\subsection {Nonequivariant homology theories.} First, we can think of $X$ 
as a nonequivariant pro-spectrum. We can apply any of the usual  
nonequivariant generalized homology functors to $X,$ such as stable 
homotopy, K-theory, or bordism. Of particular interest will be the 
ordinary homology. 

For any generalized homology $h_*,$ there is an associated dual cohomology 
theory:
$$ h^n(X) = h_{-n}(DX).$$

\subsection {Some equivariant homology theories.} The material here is 
taken from \cite{GM}. First, a bit of notation: $M_+$ is the disjoint 
union of $M$ and a basepoint, while $\tilde M$ is the unreduced suspension 
of $M$ with one of the cone points as basepoint. 

The simplest homology theory which takes into account the $S^1$ 
equivariance is \emph {Borel homology}:
$$ H^{S^1}_n(X) = H_{n-1}(ES^1_+ \wedge_{S^1} X).$$

There is also \emph {Borel cohomology}:
$$ H_{S^1}^n(X) = H^{n}(ES^1_+ \wedge_{S^1} X).$$

However, these two theories are not dual to each other, as one can easily 
see from the example of $X = \ss = D\ss,$ when both $H^{S^1}_*$ and 
$H^*_{S^1}$ are nonzero in infinitely many positive degrees but zero in 
negative degrees.

The dual homology theory to Borel cohomology is called \emph {coBorel 
homology} (or $c$-homology), and is defined by:
$$ cH^{S^1}_n(X) = \colim_m [\Sigma^{n+m} ES^1_+, K_m \wedge X]_{S^1},$$
where $K_m = K(\zz, m)$ is the Elienberg-MacLane spectrum. Notice the 
analogy with usual homology:
$$ H_n(X) = \colim_m [S^{n+m}, K_m \wedge X].$$

For example, when $X = \ss,$ one can compute $cH_n^{S^1}(\ss) = \zz$ if $n 
\leq 0$ is even and $cH_n^{S^1}(\ss) = 0$ otherwise.

Similarly, there is a dual cohomology theory to Borel homology, which is 
called \emph {coBorel cohomology} (or $f$-cohomology):
$$ fH^n_{S^1}(X) = H_{-n}^{S^1} (DX).$$

We need to introduce one more pair of dual theories: \emph {Tate homology}  
and \emph {Tate cohomology}. These were originally defined for 
spaces with finite group actions by Swan. The analogue for $S^1$-spaces 
which we use below is due to Jones \cite{J}.
$$ tH_n^{S^1}(X) = cH_n^{S^1}(\tilde ES^1 \wedge X); \hskip6pt 
tH_n^{S^1}(X) = fH^{n+1}_{S^1}(\tilde ES^1 \wedge X).$$

The group $H^*_{S^1}(\ss)= \zz[U],$ with deg $U = 2,$ acts on Borel, 
coBorel, and Tate cohomologies by cup product and on the respective Borel 
homologies by cap product. Correspondingly, the action of $U$ increases 
degree by 2 in cohomology and decreases degree by 2 in homology. We can 
also think of the ordinary homology and cohomology as $\zz[U]$ modules 
with the trivial $U$ action. 

There are long exact sequences of $\zz[U]$-modules:
$$\begin {CD}
\label {noneq}
\cdots \to H_n(X) \to H_{n+1}^{S^1}(X) @>{\operatorname{U}}>> 
H_{n-1}^{S^1}(X) \to H_{n-1}(X) \to \cdots
\end {CD}$$

\begin {equation}
\label{cbt}
\cdots \to H_n^{S^1}(X) \to cH_n^{S^1}(X) \to tH_n^{S^1}(X) \to 
H_{n-1}^{S^1}(X)\to \cdots
\end {equation}

When applied to the invariant $\swf_0,$ in light of Conjecture~\ref{osz}, 
these long exact sequences mimic the ones in Ozsv\'ath-Szab\'o theory 
from \cite{OS2}. There are also analogous sequences in cohomology.

Let us conclude with a few remarks on Tate homology and cohomology: when 
$X$ is a free $S^1$-spectrum, according to \cite{GM}:
$$ tH_*^{S^1}(X) = tH^*_{S^1}(X) =0.$$

More generally, Tate cohomology can be computed by localizing Borel 
cohomology:
$$ tH^*_{S^1}(X) = U^{-1} H^*_{S^1}(X).$$

By the localization theorem, when $X$ is a finite $S^1$-CW complex with 
semifree $S^1$ action and $X^{S^1}$ is its fixed point set, we have:
$$ tH^*_{S^1}(X) =  tH^*_{S^1}(X^{S^1}) = U^{-1} H^*_{S^1}(X^{S^1}).$$

In particular, when $Y$ is a homology 3-sphere and $\swf(Y, \spc)$ is its 
Seiberg-Witten Floer spectrum as defined in \cite{M}, for generic 
perturbations there is a unique nondegenerate reducible monopole. Thus 
$\swf^{S^1} = \ss$ and
$$ tH^*_{S^1}(X) = \zz[U, U^{-1}].$$

When $X = \swf(Y, \spc, g, A_0)$ for $b_1(Y) = 1$ and $\spc$ nontorsion as 
in Section~\ref{sec:g}, we have $X^{S^1}=*$ and
$$ tH^*_{S^1}(X) = 0.$$

\subsection {The effect of $j$ on homology.}

In Theorem~\ref{nontorsion} (b) it is claimed that $j: \swf_0 \to \swf$ 
induces isomorphisms on homology and Borel homology. Here we prove this 
claim. We let $h_*$ stand for either ordinary or Borel homology.

Recall from Subsection~\ref{sec:s0} the construction of the map $j.$ It 
comes from an attractor-repeller sequence:
$$ V(m',n) \to U(m,n) \to W(m, m'),$$
for $m \ll m' \ll 0 \ll n.$ Here $W(m, m') = B_m' \setminus A_m$ and let 
$J''(m, m')$ be the corresponding Conley index. There is a 
coexact sequence of Conley indices:
$$ J'(m',n) \to J(m,n) \to J''(m, m') \to \Sigma J'(m',n) \to \cdots $$

Applying the functor $h_*$ we get a long exact sequence:
\begin {equation}
\label {eki}
h_s(J'(m',n)) \to h_s(J(m,n)) \to h_s(J''(m, m')) \to h_{s-1}(J'(m',n)) \to 
\cdots 
\end {equation}

\begin {lemma}
\label{lemy}
Fix $k \in \zz.$ Then for every $m \ll m' \ll 0,$ we have $h_k(J''(m, 
m'))=0.$
\end {lemma}

\noindent {\bf Proof.} Because of periodicity,
$$ h_k(J''(m, m'))= h_{k+\ell}(J''(m+1, m'+1)) = \dots = h_{k-|m| \ell} 
(J''(0, m'-m)).$$

Now it suffices to show that there exists $s_0$ such that 
\begin {equation}
\label {bdd}
h_s(J''(0, m'-m)) =0 \text{ for all }s \geq s_0.
\end {equation}

For fixed $p = m'-m,$ this is true because $J''(0,p)$ is a finite 
desuspension of a finite $S^1$-CW-complex, and such complexes have their 
homology and Borel homology bounded below.  

Fix some $p_0 \gg 0$ and choose $s_0$ so that $$h_s(J''(0,p_0)) = 
h_s(J'(0,1)) =0$$
for all $s \geq s_0.$ By periodicity 
$$ h_s(J'(p, p+1)) = h_{s-p} (J'(0,1)) =0$$
for all $p \geq 0$ as well. Using the long exact sequences coming form 
attractor-repeller pairs:
$$ h_s(J'(p, p+1)) \to h_s(J'(0, p+1)) \to h_s(J'(0,p)) \to \dots $$
we obtain (\ref{bdd}) for any $p=m'-m \geq p_0$ by induction on $p.$ $\fin$

\medskip
Lemma~\ref{lemy} says that the map $h_s(J'(m',n)) \to h_s(J(m,n)) $ in 
(\ref{eki})is an isomorphism for $m \ll m' \ll 0.$ Taking direct limits as 
$n \to \infty,$ inverse limits as $m \to -\infty,$ and finally direct limits 
as $m' \to -\infty,$ we obtain that 
$$ (h_s)_* j : h_s(J'(-\infty, \infty)) \to h_s(J(-\infty, \infty))$$
is an isomorphism as well.

\subsection {The case of $S^1 \times S^2$}

Let us make concrete the difference between $\swf$ and $\swf_0$ by means of 
an example. Let $Y = S^1 \times S^2$ and $\spc$ a nontorsion $spin^c$ 
structure with $\ell > 0$ an even integer. Choose $g$ to be the product of 
the flat metric on $S^1$ and the round metric on $S^2.$ 

It is easy to see that in this case all Seiberg-Witten solutions are 
reducible. When $[\nu] = \epsilon =0,$ there are no reducibles either, so in 
fact
$$ \swf(Y, \spc, g, A_0) = *$$

Since $j$ is an isomorphism at the level of homology and Borel homology, 
for $\swf_0 = \swf_0(Y, \spc, g, A_0)$ we must have:
$$ H_*(\swf_0) = H^{S^1}_*(\swf_0) =0.$$

From the sequence (\ref{cbt}) it follows that
$$ cH^{S^1}_*(\swf_0) =tH^{S^1}_*(\swf_0).$$ 

If Conjecture~1 is true, we expect:
$$ cH^{S^1}_*(\swf_0) =tH^{S^1}_*(\swf_0) = \zz[U, U^{-1}],$$
and the periodicity map $\swf_0 \to \Sigma^{-\ell} \swf_0$ should be 
given by the action of $U^{\ell/2}.$ This computation can also be carried 
out in the context of Seiberg-Witten theory \cite{KM}.

\smallskip

\section {Relative invariants of four-manifolds with boundary}

In this section we prove Theorem~\ref{relative}. Let $X$ be a compact 
4-manifold with boundary $Y$ such that the image of $H^1(X;\rr)$ in 
$H^1(Y; \rr)$ is $0.$ Assume that $X$ has a 
$spin^c$ structure $\hat \spc$ which extends $\spc,$ and that we are given 
orientations on $H^1(X; \rr)$ and $H^2_+(X; \rr).$ 

Our goal is to construct a relative Seiberg-Witten invariant of $X$ in
the form of an element in a stable homotopy group of $\swf_0(Y, \spc).$
This construction was done in Section 9 of \cite{M} for the case
$b_1(Y)=0$, and then corrected by Khandhawit in \cite{Kh}. The current case is more or less similar, so here we will
only sketch the construction. The reader is referred to \cite{M, Kh} for the
analytical details.

Let us also choose a $spin^c$ connection $\hat A_0$ on $X$ which restricts 
to $A_0$ on $Y.$ Then we can define the Seiberg-Witten map
\begin {eqnarray*}
SW^{\mu} : i\Omega_g^1(X) \oplus \Gamma(W^+) & \to & i\Omega^2_+ (X) \oplus 
\Gamma(W^-) \oplus V^{\mu} \\
(\hat a, \hat \phi) &\to & (F^+_{\hat A_0 + \hat a}-\rho^{-1}(\phi\phi^*)_0, 
D_{\hat A_0 + \hat a}, p^{\mu} i^*(\hat a, \hat \phi))
\end {eqnarray*}
We need to explain the notation. The space $\Omega^1_g(X)$ consists of 1-forms on 
$X$ in double Coulomb gauge, as in \cite[Definition 1]{Kh}. 
Furthermore, $W^+$ and $W^-$ are the positive and negative spinor bundles on $X,$ 
respectively, $\rho$ is Clifford multiplication, and $D$ is the 
four-dimensional Dirac operator. Finally, 
$i^*$ is the restriction to $Y$, $\Pi$ is the Coulomb projection for $Y,$
and $p^{\mu}$ is the orthogonal projection to $V^{\mu} = 
V^{\mu}_{-\infty}.$ 

The map $SW$ can be decomposed into a linear and a compact map between 
suitable Sobolev completions of the domain and the target. We can apply 
Furuta's technique of finite dimensional approximation and obtain a map:
$$ SW^{\mu}_{\lambda, U} = \pr_{U \times V^{\mu}_{\lambda}} SW^{\mu} : U' 
\to U \times V^{\mu}_{\lambda}.$$
Here $U, U'$ are finite dimensional spaces. Take a small ball 
centered at the origin $B(\epsilon) \subset U$ and consider the preimage 
of $B(\epsilon) \times V^{\mu}_{\lambda}$ in $U'.$  Let $\bar K_1, \bar 
K_2$ be the intersections of this preimage with a large ball $B'$ in 
$U'$ and its boundary, respectively. Finally, map $\bar K_1, \bar K_2$ back 
to $V^{\mu}_{\lambda}$ by composing $SW^{\mu}_{\lambda, U}$ with the 
obvious projection. Denote the respective images by $K_1$ and $K_2.$

Let us assume for the moment that we are in the simplest case, when $b_1(Y) 
=1$ and $c_1(\spc)$ is not torsion. Recall the notations from Section 3. 
Since 
$K_1$ is compact, we can choose $n, -m,$ and $R$ sufficiently large so that 
$K_1 \subset T(R).$ Furthermore, the analysis done in \cite{M} (based on 
the compactness properties of the four-dimensional Seiberg-Witten 
equations) shows that there exists an index pair $(N, L)$ for the invariant 
part of $T(2R)$ in the gradient flow such that $K_1 \subset L$ and $K_2 
\subset N.$

Thus we can define a map
$$ (U')^+ \cong B'/\del B' \to (B(\epsilon) \times N)/(B(\epsilon) \times L 
\cup \del B(\epsilon) \times N) \cong U^+ \wedge I^{\mu}_{\lambda}$$
by applying $SW^{\mu}_{\lambda, U}$ to the elements of $\bar K_1$ and 
sending everything else to the basepoint. 

For $-\lambda, \mu \gg 0,$ after taking the relevant desuspensions, this 
gives an element in a stable homotopy group of $J = J(m,n).$ There are such 
elements for any $-m, n \gg 0,$ and they commute with the maps between the 
different $J(m,n)$ coming from the attractor-repeller exact triangles. 
Therefore, we can take direct limits and pro-limits and obtain a map to 
$\swf(Y, \spc).$ If we insist on doing the constructions equivariantly, we 
get an element
$$ \Psi(X, \hat \spc, \hat A_0) \in \pi_*^{S^1}(\swf(Y, \spc)). $$ 

In the end we find that for any $X$ and $Y$ there is an element:
$$  \Psi(X, \hat \spc, \hat A_0) \in \pi_*^{S^1}(\swf(Y, \spc)). $$

One can show by a continuity argument that $\Psi$ is independent of the 
choices made in its construction, up to canonical equivalence. 

Starting from here we can compose with the canonical map from stable 
homotopy to any other generalized homology theory $h.$ Thus we obtain 
relative Seiberg-Witten invariants of $X$ with values in any
$h_*(\swf(Y,\spc)).$

Now let us vary the base connection $\hat A_0$ 
on $X$ by adding to it a harmonic 1-form $\alpha\in H^1(X; \rr)$ which 
annihilates the normal vector to the boundary. 

We can collect together the maps $\Psi(X, \hat\spc, \hat A_0+\alpha)$ 
as $\alpha$ varies over the Picard torus $P=H^1(X; \rr) / H^1(X; \zz)$ and 
obtain a bundle of morphisms from spheres 
to $\swf.$ In other words, we get a morphism from the Thom space of a 
vector bundle over $P$ (the Dirac index bundle) to $\swf(Y, \spc):$
$$ \Psi(X, \hat\spc) : T(Ind) \to \swf(Y, \spc).$$

In the case when $X$ is closed, this is the invariant constructed by 
Bauer in \cite{B}. 

\begin {remark} 
By adding a nonexact 2-form perturbation together with a small $\epsilon 
\sin$ perturbation to the Seiberg-Witten map, an analogous argument should 
give a morphism from the Picard torus to $\swf_0.$ Thus, we expect $\Psi$
to factor through the map $j: \swf_0 \to \swf.$
\end {remark}

\appendix 
\section {The K-theoretic obstruction.}

Let $Y$ be a closed, oriented, Riemannian 3-manifold, endowed with a
$spin^c$ structure $\spc.$ We keep the notations from Section~2, but do 
not impose any condition on $b_1$ or on $\spc.$

In particular, if we fix a base connection $A_0,$ we identify all other 
connections $A$ with forms $a \in i\Omega^1(Y; \rr)$ via $a = A - A_0.$
For each $a$ there is a Dirac operator $\dir_a.$ We restrict our attention 
to harmonic 1-forms $a \in H^1(Y; \rr).$ In particular, if $a \in H^1(Y; 
\zz),$ there is a harmonic map $u:Y \to S^1$ with $a=u^{-1}du.$ In this 
case, for every spinor $\phi \in \Gamma(W),$ 
$$ \dir_a \phi = \dir_0 (u\phi).$$

We can then form a Hilbert bundle over the Picard torus
$$P = H^1(Y; \rr)/ H^1(Y; \zz),$$
with fiber $\Gamma(W)$ by gluing via isomorphisms of the form $\phi \to 
u\phi.$ In this setting we get a continuous family of Dirac operators 
$\dir_a$ acting on each fiber, which is parametrized by $P.$ 

As explained in \cite{PS} and \cite{CJS}, such a family induces a family 
of polarizations over $P,$ and therefore a \emph {structural map} $P \to 
U,$ where $U$ denotes the infinite unitary group. 

Denote:
$$ H^*(Y,\spc) = c_1(\spc)^{\perp} = \{a \in H^*(Y): a \cup
c_1(\spc)=0\};$$
$$\tilde P = H^1(Y, \spc; \rr)/ H^1(Y, \spc; \zz).$$
Let $$b_1(Y,\spc) = \dim H^1(Y,\spc; \rr).$$ 
Note that $b_1(Y, \spc)$
equals $b_1(Y)$ when $c_1(\spc)$ is torsion and $b_1(Y) -1$ otherwise.

The obvious covering map $\tilde P \to P$ corresponds to the group 
homomorphism $\pi_1(P) \to \pi_1(U) = \zz$ given by the spectral flow.

We lift the family of Dirac operators from $P$ to $\tilde P.$ The 
composition $\tilde P \to P \to U$ induces an element 
$$q(Y, \spc) \in K^1(\tilde P).$$ 
According to \cite{CJS}, this is the 
obstruction for the Seiberg-Witten flow category to be framed, or, in 
other words, for Floer stable homtopy to be well-defined. 

Note that since $\tilde P$ is a torus, by the K\"unneth
formula 
$$ K^1(\tilde P) = H^1(\tilde P) \oplus H^3(\tilde P) \oplus 
H^5(\tilde P) \oplus \dots.$$

\begin {proposition}
The obstruction $q(Y, \spc) = q \in K^1(\tilde P)$ is
given by the intersection form:
$$ \Lambda^3 H^1(Y, \spc; \rr) \to \rr, \hskip6pt (a_1, a_2, a_3) \to (a_1
\cup a_2 \cup a_3)[Y],$$
considered as an element in $(\Lambda^{3}H^1(Y, \spc; \rr))^* = 
H^3(\tilde P) \subset K^1(\tilde P).$
\end {proposition}

\noindent{\bf Proof.} According to \cite{APS}, the obstruction $q \in 
K^1(\tilde P) = K^0(\tilde P \times S^1)$ is given by the K-theoretic 
index 
of the family $\{D_{t,a}\}$ of Dirac operators on $Y \times S^1:$ 
$$ D_{t,a} = \begin{cases}
I\cos t  + i(\dir_a + \delt ) \sin t, &  0 \leq t \leq \pi;\\
(\cos t + i\sin t)I, &  \pi \leq t \leq 2\pi.
\end{cases}
$$
parametrized by $(a, t) \in \tilde P \times S^1.$ 

Since $K^0(\tilde P \times S^1)$ is nontorsion, it suffices to compute
the Chern character of $q \in H^*(\tilde P \times S^1).$ This can be 
done
using the index theorem for families (\cite{AS}). Set $m= b_1(Y,
\spc).$ Let $a_1, \dots, a_{m}$ be a basis for $H^1(Y, \spc; \zz),$
$a_{m+1}$ a generator of $H^1(S^1; \zz),$ and
$\alpha_1, \dots, \alpha_{m+1}$ the dual basis for $H^1(\tilde P \times 
S^1).$
The Dirac line bundle $\mathcal{L}$ over $Y\times S^1 \times \tilde P 
\times S^1$
has first Chern class:
$$ c_1(\mathcal{L}) = \sum_{i=1}^{m+1} \alpha_i a_i.$$

The index theorem gives:
$$ ch(q) = e^{\sum \alpha_i a_i} \hat{A}(Y
\times S^1) [Y \times S^1].$$

The 3-manifold $Y$ is parallelizable, hence so is $Y \times S^1.$ Thus the
$\hat A$ genus is $1.$

In the Taylor decomposition the only term that survives is:
$$ ch(q) = \frac{1}{4!}(\sum_{i=1}^{m+1} \alpha_i a_i)^4 [Y \times S^1],$$
which corresponds to the intersection form in $H^4(Y \times S^1) \cong
H^3(Y).$

\medskip
\noindent\textbf{Acknowledgements. } The second author would like to thank 
Daniel Biss, Danny Calegari, Mikio Furuta, and Michael Hopkins
for helpful discussions during the preparation of this paper. 
\smallskip

\end {document}